\documentclass{article}
\usepackage{latexsym}
\usepackage{graphicx}
\usepackage{amsmath}
\usepackage{amsfonts}
\usepackage[latin1]{inputenc}
\usepackage{amssymb}
\usepackage{color}

\newtheorem{theo}{Theorem}
\newtheorem{coro}{Corollary}
\newtheorem{ex}{Example}
\newtheorem{lemm}{Lemma}
\newtheorem{prop}{Proposition}
\newtheorem{definition}{Definition}

\title{Bivariate copulas defined from matrices}
\author{C\'{e}cile Amblard$^{(1)}$, St\'ephane Girard$^{(2,\star)}$, Ludovic Menneteau$^{(3)}$}

\date{\small $^{(1)}$ Universit\'e Grenoble 1, {\tt Cecile.Amblard@imag.fr}\\
 $^{(2,\star)}$ Inria Grenoble Rh\^one-Alpes \& LJK, {\tt Stephane.Girard@inria.fr}\\
 $^{(3)}$ Universit\'e Montpellier 2, {\tt mennet@math.univ-montp2.fr}}
 
\begin{document}

\maketitle

\begin{abstract}
We propose a semiparametric family of copulas based on a set of orthonormal functions
and a matrix. This new copula permits to reach values of Spearman's Rho
arbitrarily close to one without introducing a singular component. 
Moreover, it  encompasses several extensions of FGM copulas as well 
as copulas based on partition of unity such as Bernstein or checkerboard copulas.
Finally, it is also shown that projection of arbitrary densities of copulas
onto tensor product bases can enter our framework.\\



\end{abstract}

\section{Introduction}
%---------------------------------------------------------------%

A bivariate copula defined on the unit square $[0,1]^2$ is a bivariate cumulative distribution function (cdf)
with univariate uniform margins.
Sklar's Theorem~\cite{SKLAR} states that any
bivariate distribution with cdf $H$
and marginal cdf $F$ and $G$ can be written
$H(x,y)=C(F(x),G(y))$, where $C$ is a copula. 
This result justifies the use of copulas for building bivariate distributions.
One of the most popular parametric family of copulas is the Farlie-Gumbel-Morgen\-stern
(FGM) family~\cite{FAR,GUMBEL,MORG} defined when $\theta\in[-1,1]$ by
\begin{equation}
\label{eqFGM}
C(u,v)=uv+\theta u(1-u)v(1-v).
\end{equation}
A well-known limitation to this
family is that it does not allow the modeling of large dependences
since the associated Spearman's Rho is limited to $[-1/3,1/3]$.
A possible extension of the FGM family is to consider
the semi-parametric family of symmetric copulas defined by
\begin{equation}
\label{eqcopunous}
C(u,v)=uv+\theta \varphi(u)\varphi(v),
\end{equation}
with $\theta\in[-1,1]$. It was first introduced in~\cite{LLALENA},
and extensively studied in~\cite{Nous,Nous2}. In particular, it
can be shown that, for a properly chosen function $\varphi$,
the range of Spearman's Rho is extended to $[-3/4,3/4]$.
In~\cite{Metrika} an extension of~(\ref{eqcopunous}) is proposed
where $\theta$ is a univariate function. This modification
allows the introduction of a singular component concentrated on
the diagonal $v=u$ and extends the range of Spearman's Rho to $[-3/4,1]$.
We also refer to~\cite{DUR,DUR2} for another yet similar extensions.
\\

\noindent Here, a new extension of~(\ref{eqcopunous}) is proposed
where, roughly speaking, the single parameter $\theta$ is replaced
by a matrix and the function $\varphi$ is replaced by a set of functions.
This new copula permits to reach values of Spearman's Rho
arbitrarily close to 1 without singular component. 
Moreover, it also encompasses copulas based on partition of unity such
as Bernstein copula~\cite{sancetta} or checkerboard copula~\cite{LMST,LMT}.
Finally, it is also shown that projection of arbitrary densities of copulas
onto tensor product bases can enter our framework.
This paper is organized as follows: The family of copula is introduced
in Section~\ref{secdef} and some first properties are given. Algebraic properties
are established in Section~\ref{secalge}, dependence properties are reviewed in
Section~\ref{secdep} while important approximation issues are highlighted in 
Section~\ref{secapp}. Some links with existing copulas as well as examples of new copulas are 
presented in Section~\ref{secex}. Concluding remarks are drawn Section~\ref{seconcl}.
Proofs are postponed to the Appendix.

\section{A new family of copulas}
%---------------------------------------------------------------%
\label{secdef}

Throughout this paper, $e_j$ denotes the $j$th vector of the canonical
basis of ${\mathbb R}^p$ where $j=1,\dots, p$ and $p\geq 2$. Besides, $\langle f,g \rangle=\int_0^1 f(t) g(t) dt$ is the usual scalar product in $L_2([0,1])$. The associated norm is denoted by $\|f\|=\langle f,f \rangle^{1/2}$.
We focus on the modeling of the copula density (denoted by $c$) rather than on the modeling of the copula (denoted by $C$) itself:
\begin{definition}
\label{def1}
Let $\phi:[0,1]\to {\mathbb R}^p$ be a vector of $p$ orthonormal functions such 
that $\phi_1(t)=1$ for all $t\in[0,1]$. Two sets are defined from $\phi$:
\begin{eqnarray*}
{\mathcal A}_\phi& = & \left\{A\in{\mathbb R}^{p\times p},\;A e_1=e_1, \;^t\!Ae_1=e_1,
 \forall (u,v)\in [0,1]^2,\;^t\phi(u)A \phi(v)\geq 0\right\},\\
{\mathcal C}_\phi& = &\left\{c: [0,1]^2\rightarrow {\mathbb R}, \;
c(u,v)= \;^t\phi(u) A \phi(v),\, A\in {\mathcal A}_\phi\right\}.
\end{eqnarray*}
\end{definition}
The next result establishes that all the functions of ${\mathcal C}_\phi$ are densities of
copulas.
\begin{prop}
\label{propdens}
${\mathcal C}_\phi$ is a non-empty set of copula densities.
\end{prop}
It is clear that ${\mathcal A}_\phi$ is not empty since $A_1=\;^t e_1 e_1 \in {\mathcal A}_\phi$.
The associated function $c(u,v)= \;^t\phi(u) A_1 \phi(v) = 1 \in {\mathcal C}_\phi$ is the density of the independent
copula. The remainder of the proof is postponed to the Appendix.
Let us also note 
${\mathcal A}_\phi$ is a subset of matrices with eigenvalue 1 associated to the eigenvector $e_1$
and with non-negative trace.
Indeed, if $A\in{\mathcal A}_\phi$ then
$$
\mbox{tr}(A \phi(x)\; ^t\phi(x))= \mbox{tr}(^t\phi(x)A \phi(x)) = \; ^t\!\phi(x)A \phi(x) \geq 0
$$
for all $x\in[0,1]$. Integrating with respect to $x$ yields the result since, by assumption,
\begin{equation}
\label{ortho}
 \int_0^1 \phi(x)\; ^t\phi(x) dx= I_p,
\end{equation}
where $I_p$ is the $p\times p$ identity matrix.

\begin{ex} 
\label{ex1}
 If $p=2$ then $A\in {\mathcal A}_\phi$ implies that $A=$diag$\{1,\theta\}$ with $\theta\geq -1$.
 The associated density of copulas can be written as $c(u,v)=1+\theta\phi(u)\phi(v)$
 which corresponds to family~(\ref{eqcopunous}). 
\end{ex}
This family includes FGM copulas~(\ref{eqFGM})
which contains all copulas with both horizontal and vertical quadratic 
sections~\cite{QUESADA}, the subfamily 
of symmetric copulas with cubic sections proposed in~\cite{NELSEN97},
equation~(4.4), and some kernel extensions of FGM copulas introduced in~\cite{HUE,LAI}.
We refer to~\cite{FISC} for a method to
construct admissible functions $\phi$. \\

\noindent The following lemma will reveal useful to build densities of copulas in ${\mathcal C}_\phi$
without the orthogonality assumption on $\phi$.
\begin{lemm}
\label{lemortho}
 Let $\psi:[0,1]\to {\mathbb R}^p$ be a vector of $p$ functions such 
that $\psi_1(t)=1$ for all $t\in[0,1]$ and $\int_0^1\psi(x)dx=e_1$. 
Let $\Gamma$ be the Gram matrix defined as $\Gamma=\int_0^1 \psi(x)\; ^t\!\psi(x)dx$
and $B\in {\mathcal A}_\psi$. Then, $A:=\Gamma^{1/2} B \Gamma^{1/2} \in {\mathcal A}_\phi$
where $\phi:=\Gamma^{-1/2} \psi$ fulfills the conditions of Definition~\ref{def1}
and $^t\!\phi(u)A\phi(v)=\;^t\!\psi(u)B\psi(v)$ for all $(u,v)\in[0,1]^2$.
\end{lemm}
See Subsection~\ref{proofaux} of the Appendix for a proof. A direct application of this lemma yields:
\begin{ex}
\label{excub}
The parametric family 
of copulas with cubic sections proposed in~\cite{NELSEN97}, Theorem~4.1
and given by
$$
C(u,v)=uv + uv(1-u)(1-v)[ A_1v(1-u)+A_2(1-v)(1-u) + B_1uv + B_2u(1-v)]
$$
can be written in our formalism with Lemma~\ref{lemortho}. Here, $p=3$, $\psi_1(t)=1$, $\psi_2(t)=1-4t+3t^3$,
$\psi_3(t)=2t-3t^2$, 
$\Gamma= \begin{pmatrix} 1 &0 & 0 \\ 0 & 2/15 & 1/5  \\ 0 &1/5 & 4/5 \end{pmatrix}$
and 
$B=\begin{pmatrix} 1 &0 & 0 \\ 0 & A_2 & A_1  \\ 0 &B_2& B_1 \end{pmatrix}$.
\end{ex}
More generally, iterated FGM families~\cite{IFGM1,IFGM2,IFGM3} where
  $$
C(u,v)=uv + \sum_{j=1}^p \theta_j (uv)^{\alpha_j} ((1-u)(1-v))^{\beta_j},
$$
and $\{\alpha_j,\beta_j\}=\{[j/2]+1,[(j+1)/2]\}$ can be shown to be particular cases
of our family thanks to Lemma~\ref{lemortho}.
More examples are given in Section~\ref{secex}. \\

\noindent Let $T_\phi$ the mapping defined for any copula density $c\in L_2([0,1]^2)$ by
\begin{equation}
\label{deftphi}
 T_\phi(c)=\int_0^1\int_0^1 c(x,y) \phi(x)\;^t\!\phi(y)dxdy \in {\mathbb R}^{p\times p}.
\end{equation}
The mapping $T_\phi$ permits to compute the matrix associated to any copula density $c\in {\mathcal C}_\phi$:
\begin{prop}
\label{propT}
Each copula density $c\in {\mathcal C}_\phi$ is defined by an unique matrix $A$ which is given
by $$
A=T_\phi(c)={\mathbb E}_c(\phi(U)\;^t\!\phi(V))=cov(\phi(U),\phi(V))+ e_1\; ^t\! e_1,
$$
where $(U,V)$ is a random pair with density $c$.
\end{prop}

\section{Algebraic properties}
%---------------------------------------------------------------%
\label{secalge}

Let $\times$ denote the matrix product and let $\star$ denote the product of copulas introduced in~\cite{Darsow} and defined in terms of densities as 
$$
c_A\star c_B(u,v):=\int_0^1 \; c_A(u,s)c_B(s,v)ds,
$$
for all $(u,v)\in[0,1]^2$. The following stability properties can be established (see Appendix for a proof).
\begin{prop}
\label{propA}
 ${\mathcal A}_\phi$ is a convex set and $({\mathcal A}_\phi,\times)$ is a semi-group.
 If, moreover, $I_p \in {\mathcal A}_\phi$ then $({\mathcal A}_\phi,\times)$ is a monoid.
\end{prop}
Let us also consider the bivariate function 
$$
q(u,v)= \;^t\!\phi(u)\phi(v)
$$
defined for $(u,v)\in[0,1]^2$.
The following result is the analogous of Proposition~\ref{propA}
for ${\mathcal C}_\phi$.
\begin{prop}
\label{propC}
 ${\mathcal C}_\phi$ is a closed convex set and $({\mathcal C}_\phi,\times)$ is a semi-group.
 If, moreover, $q \in {\mathcal C}_\phi$ then $({\mathcal C}_\phi,\times)$ is a monoid.
\end{prop}
In view of Propositions~\ref{propT} --~\ref{propC}, it appears that $T_\phi(c_A\star c_B)=T_\phi(c_A) T_\phi(c_B)$ for all $(c_A,c_B) \in {\mathcal C}_\phi^2$ and thus:
\begin{prop}
$T_\phi$ is an isomorphism between $({\mathcal A}_\phi,\times)$ and $({\mathcal C}_\phi,\star)$.
\end{prop}
To summarize, it appears that ${\mathcal A}_\phi$ is stable with respect to
matrix multiplication. Moreover, multiplying the matrices 
is equivalent to ``multiplying'' the copulas using the $\star$ product.\\

\noindent Finally, the next lemma shows that it is possible to aggregate copulas
of ${\mathcal C}_\phi$ with different number of orthogonal functions
through the use of Ces\`aro summations.
\begin{lemm}
 \label{lemagg}
Let $c_p\in {\mathcal C}_\phi$ with associated matrix $A$ and consider the 
density of copula defined for all $(u,v)\in[0,1]^2$ and $q\geq 1$ by
$$
\bar{c}_{q}(u,v):= \frac{1}{q} \sum_{p=1}^{q} c_{p}(u,v).
$$
Then, $\bar{c}_{q}\in {\mathcal C}_\phi$ with associated matrix $B$ defined by
$B_{ij}= (q+1 - \max(i,j)) A_{i,j}/q$ for all $(i,j)\in\{1,\dots,q\}^2$.
\end{lemm}
An application of this lemma is illustrated in Paragraph~\ref{trigo}.

\section{Dependence properties}
%---------------------------------------------------------------%
\label{secdep}

Several measures of association between the components of a random pair 
can be considered:
the normalized volume~\cite{SIGMA},
Kendall's Tau~\cite{NELSEN06}, paragraph~5.1.1,
Gini's gamma~\cite{JNS10}, Blomqwist's medial
correlation coefficient~\cite{NELSEN06}, paragraph~5.1.4,
Spearman's footrule~\cite{JNS10},
and Spearman's Rho~\cite{NELSEN06}, paragraph~5.1.2.
All these measures are
invariant to strictly increasing functions.
Kendall's Tau and Spearman's Rho can be interpreted as probabilities of concordance minus 
probabilities of discordance of two random pairs.
Let us first focus on the Spearman's Rho. 
It can be written only in terms of the copula $C$:
\begin{equation}
 \label{eqrho}
 \rho=12\int_0^1 \!\!\! \int_0^1 C(u,v) dudv\!-\!3.
\end{equation}
Note that $\rho$ coincides with the correlation coefficient between the uniform
marginal distributions.
In the case of the copula introduced in Definition~\ref{def1}, 
it can be expressed thanks to the function $\phi$ and the 
matrix $A\in {\mathcal A}_\phi$ associated to its density.
\begin{prop}
\label{proprho}
Let $(U,V)$ be a random pair with density of copula $c\in{\mathcal C}_\phi$
associated to the matrix $A\in{\mathcal A}_\phi$.
The Spearman's Rho is given by $\rho=12\;^t\!\mu A\mu-3$
where $\mu =\int_0^1x\phi(x)dx$.
\end{prop}
Similarly to the Spearman's Rho, the Kendall's Tau can be written only in terms of the copula $C$:
\begin{equation}
\label{eqtau}
\tau=4\int_0^1 \!\!\! \int_0^1 C(u,v) dC(u,v)\!-\!1,\;
\end{equation}
and in the framework of Definition~\ref{def1}, 
it can be expressed thanks to the function $\phi$ and the 
matrix $A\in {\mathcal A}_\phi$ associated to its density.
\begin{prop}
\label{proptau}
Let $(U,V)$ be a random pair with density of copula $c\in{\mathcal C}_\phi$
associated to the matrix $A\in{\mathcal A}_\phi$.
The Kendall's Tau is given by $\tau=1-4 \mbox{tr} \left( \;^t\!A \Theta A \Theta\right)$
where $\Theta$ is the $p\times p$ matrix defined by $\Theta=\int_0^1 \Psi(u)\;^t\!\phi(u)du$
and with $\Psi(u)=\int_0^u \phi(t)dt$.
\end{prop}
Let us note that Propositions~\ref{proprho} and~\ref{proptau} extend the results
of~\cite{GRO}, Theorem~25 established in the case of copulas based on partition of unity.
Such copulas were introduced in~\cite{LMST,LMT}. It is shown in Section~\ref{secex}
that they are particular cases of the family considered in Definition~\ref{def1}.
Besides, Blomqwist's medial correlation coefficient, Gini's gamma and Spearman's footrule 
can also be rewritten in terms of the copula $C$ (see for instance~\cite{JNS98}). All these coefficients
thus benefit from closed form expressions similar to these of Propositions~\ref{proprho},~\ref{proptau}
and we do not enter into details there.\\

\noindent The upper tail dependence between two random variables $X$ and $Y$ with respective
cdf $F$ and $G$ can be quantified
by the upper tail dependence coefficient~\cite{JOE}, paragraph~2.1.10,  defined as: 
$$
\lambda=\lim_{t\to 1^-} {\mathbb P}(F(X)>t | G(Y)>t).
$$
Again, this coefficient can be written only in terms of the copula 
(see~\cite{NELSEN06}, Theorem~5.4.2):
\begin{equation}
 \label{eqlambda}
 \lambda=\lim_{u\rightarrow 1^-}\frac{\bar{C}(u,u)}{1-u},
\end{equation}
where $\bar{C}$ is the survival copula, {\it i.e.} $\bar{C}(u,v)=1-u-v+C(u,v)$.
In the family ${\mathcal C}_\phi$, the upper tail dependence is not possible:
\begin{prop}
\label{proplambda}
Let $(U,V)$ be a random pair with density of copula $c\in{\mathcal C}_\phi$.
Then, $\lambda=0$.
\end{prop}

\section{Approximation properties}
%---------------------------------------------------------------%
\label{secapp}

Let $c$ be a density of copula in $L_2([0,1]^2)$. Recall that the mapping $T_\phi$ associates
a $p\times p$ matrix to $c$ via~(\ref{deftphi}) and introduce:
\begin{equation}
\label{defP}
  P_\phi(c)(u,v)=\;^t\!\phi(u) T_\phi(c)\phi(v), \;\;(u,v)\in[0,1]^2.
\end{equation}
The next lemma gives a necessary and sufficient condition for $P_\phi(c)\in {\mathcal C}_\phi$:
\begin{lemm}
\label{lemP}
\hspace*{\fill}
\begin{itemize}
 \item [(i)] Let $c\in L_2([0,1]^2)$ be an arbitrary density of copula.
If $I_p\in {\mathcal A}_\phi$ then $P_\phi(c)\in {\mathcal C}_\phi$ 
and $P_\phi(c)=q\star c \star q$.
\item [(ii)] Conversely, if $P_\phi(c)\in {\mathcal C}_\phi$ for all $c\in L_2([0,1]^2)$
then $I_p\in {\mathcal A}_\phi$.
\end{itemize}
\end{lemm}
Let $(c_1,c_2)\in L^2_2([0,1]^2)$ and let us consider the scalar product defined as
$$
\prec c_1,c_2\succ = \int_0^1 (c_1\star c_2)(u,u)du = \int_0^1\int_0^1 c_1(u,v) c_2(v,u) du dv.
$$
Let us remark that, for symmetric copulas, the above scalar product reduces to the $L_2-$ scalar product
$$
\langle c_1,c_2\rangle = \int_0^1\int_0^1 c_1(u,v) c_2(u,v) du dv.
$$
In the case of densities of copulas in ${\mathcal C}_\phi$, 
the scalar product can be computed using the associated matrices:

\newpage

\begin{lemm}
\label{lemps}
 Let $c_2\in L_2([0,1]^2)$.
 \begin{itemize}
  \item [(i)] If $c_1\in {\mathcal C}_\phi$ with associated matrix $A$, then 
  $\prec c_1,c_2\succ = \mbox{tr}( AT_\varphi(c_2) )$.
  \item [(ii)] If, moreover, $c_2\in {\mathcal C}_\phi$ with associated matrix $B$, then 
  $\prec c_1,c_2\succ = \mbox{tr}( AB )$.
 \end{itemize}
\end{lemm}
As a consequence of the above lemmas, we have:
\begin{prop}
 \label{propP}
$P_\phi$ is an orthogonal projection on ${\mathcal C}_\phi$ if and only if $I_p\in {\mathcal A}_\phi$.
\end{prop}
This result is now illustrated on the FGM family where explicit computations can be done:
\begin{ex}
It is well known that the FGM family
is a particular case of Example~\ref{ex1}
with $\phi(x)=\sqrt{3}(1-2x)$, $A=\mbox{diag}\{1,\theta\}$ and where $|\theta |\leq1/3$.
Here, $I_p\notin {\mathcal A}_\phi$ and thus 
the projection $P_\phi(c)$ of any density of copula $c$ on the FGM family 
is not itself a density of copula in the general case. However, let us remark that
$T_\phi(c)=\mbox{diag}\{1,\widetilde{\theta}\}$
and thus 
$P_\phi(c)(u,v)=1+3\tilde\theta(1-2u)(1-2v)$
where 
\begin{eqnarray*}
\tilde{\theta}=\int_0^1\!\!\int_0^1 c(x,y)\phi(x)\phi(y)dxdy =3\int_0^1\!\!\int_0^1 c(x,y)(1-2x)(1-2y)dxdy =\rho_c,
\end{eqnarray*}
the Spearman's Rho associated to $c$.
We thus have the following result: 
\begin{itemize}
\item If $|\rho_c|\leq 1/3$ then $P_\phi(c)$ is a FGM copula and $\rho_{P_\phi(c)}=\rho_c$,
\item If $|\rho_c|>1/3$ then $P_\phi(c)$ is not a copula.
\end{itemize}
\end{ex}
It appears from this example that it is possible to associate to any copula 
a FGM copula with the same Spearman's Rho $\rho_c$ provided $|\rho_c|\leq 1/3$.\\

\noindent Suppose now that $\{\phi_i\}_{i\geq 1}$ is an orthonormal basis of $L_2([0,1])$.
Then, $\{\phi_i \otimes \phi_j \}_{i,j\geq 1}$ is an orthonormal basis of $L_2([0,1]^2)$
where $\otimes$ denotes the tensor product, {\it i.e} $(f\otimes g)(u,v):=f(u) g(v)$ for all $(u,v)\in [0,1]^2$. As a consequence,
the $L_2-$projection of any $c\in L_2([0,1]^2)$ on $\{\phi_i \otimes \phi_j \}_{1\leq i,j\leq p}$
is given by
$$
\tilde{c}_p(u,v)=\sum_{i=1}^p \sum_{j=1}^p a_{i,j} \phi_i(u) \phi_j(v) = \;^t\!\phi(u) A \phi(v)
$$
where $A=(a_{i,j})_{1\leq i,j\leq p}$ with 
$$
a_{i,j}=\int_0^1\!\!\int_0^1 c(x,y) \phi_i(x) \phi_j(y) dxdy = (T_\phi(c))_{i,j}.
$$
This yields $A=T_\phi(c)$ and $\tilde{c}_p=P_\phi(c)$. In view of Lemma~\ref{lemP}(i), 
it follows that $\tilde{c}_p\in {\mathcal C}_\phi$ if $I_p\in {\mathcal A}_\phi$.
As a conclusion, the $L_2$-projection of any density of copula in $L_2([0,1]^2)$ on a tensor product basis
can be written in our formalism and the following result holds:

\newpage

\begin{theo}
\label{propproj}
Let $\{\phi_i\}_{1\leq i\leq p}$ be an orthonormal family of $L_2([0,1])$.
\begin{itemize}
 \item [(i)] The projection $P_\phi$ on ${\mathcal C}_\phi$ introduced in~(\ref{defP})
coincides with the $L_2-$ projection on $\{\phi_i \otimes \phi_j \}_{1\leq i,j\leq p}$. 
\item [(ii)] Moreover, these projections give rise to densities of copula in ${\mathcal C}_\phi$ if and only if
$I_p\in {\mathcal A}_\phi$.
\end{itemize}
\end{theo}
Classical approximation properties in $L_2([0,1]^2)$ yield:
\begin{coro}
 \label{coroprop}
Suppose the assumptions of Theorem~\ref{propproj} hold and let $c\in L_2([0,1]^2)$.\\
Then, $\|c-P_\phi(c)\|\to 0$ and $\rho(P_\phi(c))\to \rho(c)$ as $p\to\infty$, where 
$\rho(P_\phi(c))$ and $\rho(c)$ denote respectively the Spearman's Rho associated to $P_\phi(c)$ and $c$.
\end{coro}

\section{Examples}
%---------------------------------------------------------------%
\label{secex}

In paragraph~\ref{expart1}, some examples of copulas found in the literature 
are shown to enter our model. New families are exhibited in paragraph~\ref{expart2}.

\subsection{Copulas based on partition of unity}
\label{expart1}

Recall that a collection of functions $\xi=\;^t\!(\xi_1,\dots,\xi_p)$ is called
a partition of unity~\cite{LMST,LMT} if
\begin{itemize}
 \item $\xi_i\geq 0$ for all $i=1,\dots,p$,
 \item $\int_0^1 \xi_i(x) dx = 1/p$ for all $i=1,\dots,p$ and
 \item $\sum_{i=1}^p \xi_i = 1$.
\end{itemize}
It can be established that copulas based on partition of unity are particular of 
the proposed family:
\begin{prop}
 \label{propparti}
 Let $\xi=\;^t\!(\xi_1,\dots,\xi_p)$ be a partition of unity,  and let $M$ be a $p\times p$
 doubly stochastic matrix. Then, the function defined for all $(u,v)\in[0,1]^2$ by
 $c(u,v) = p \;^t\!\xi(u) M \xi(v)$ is a density of copula and $c\in {\mathcal C}_\phi$.
 Moreover, 
 $\phi= (H\Gamma_\xi \;^t H)^{-1/2} H\xi$ where $s=e_1+\dots+e_p$, $H=I_p + e_1\;^t\!s - s\;^t\!e_1$,
 and $\Gamma_\xi$ is the Gram matrix associated to $\xi$.
\end{prop}
As an illustration, we have:
\begin{ex} 
\label{exbern}
The Bernstein copula~\cite{sancetta} is obtained by choosing 
$$
\xi_i(x)= C_{p-1}^{i-1} x^{i-1} (1-x)^{p-i}
$$ 
and 
\begin{equation}
 \label{eqbern}
M_{ij} = p \left\{C\left(\frac{i}{p},\frac{j}{p}\right) -  C\left(\frac{i-1}{p},\frac{j}{p}\right)
- C\left(\frac{i}{p},\frac{j-1}{p}\right) + C\left(\frac{i-1}{p},\frac{j-1}{p}\right) \right\},
\end{equation}
where $C$ is an arbitrary copula.
\end{ex}
In the case of Bernstein copula, the basis $\phi= (H\Gamma_\xi \;^t H)^{-1/2} H\xi$ cannot be simplified.
However, in the particular case where $\{\xi_1,\dots,\xi_p\}$ is orthogonal and $\int_0^1\xi_i^2(t)dt=\beta^2$
for all $i=1,\dots,p$, we have $\Gamma_\xi=\beta^2 I_p$ and therefore 
$$
H\Gamma_\xi \;^t H =\beta^2 \begin{pmatrix}  p      & 0       & \dots & \dots & \dots & 0\\
					     0      & 2       &  1    & \dots & \dots & 1\\
					     \vdots & 1       &  2    & \ddots&       & \vdots\\
					     \vdots & \vdots  & \ddots& \ddots& \ddots& \vdots \\
					     \vdots & \vdots  &       & \ddots&  2    & 1 \\
					     0      & 1       & \dots & \dots &  1    & 2 
\end{pmatrix}.
$$
The inverse square-root of this matrix benefits from a closed-form expression,
and we thus have a simple linear relation between the two families of functions: $\phi= \Omega \xi$ with
\begin{equation}
 \label{eqomega}
\Omega=\frac{1}{\beta}      \begin{pmatrix}  p^{-1/2}      & \dots   & \dots & \dots & \dots & p^{-1/2}\\
					     -p^{-1/2}      & \gamma  & (\gamma-1)    & \dots & \dots & (\gamma-1)\\
					     \vdots & (\gamma-1)       & \gamma& \ddots&       & \vdots\\
					     \vdots & \vdots  & \ddots& \ddots& \ddots& \vdots \\
					     \vdots & \vdots  &       & \ddots&\gamma & (\gamma-1) \\
					     -p^{-1/2}      & (\gamma-1)       & \dots & \dots &  (\gamma-1)    & \gamma
\end{pmatrix},
\end{equation}
where $\gamma=(p-2+p^{-1/2})/(p-1)$. Explicit computations can be achieved in case of the checkerboard copula:
\begin{ex} 
\label{excheck}
The checkerboard copula~\cite{LMST,LMT} is obtained by choosing 
$
\xi_i(x)= {\mathbb{I}}\left\{x \in I_i\right\}
$
where $\{I_i,\; i=1,\dots,p\}$ is the equidistant partition of $[0,1]$ into $p$ intervals
$$
I_i=\left[\frac{i-1}{p},\frac{i}{p}\right]
$$
and 
$M_{ij}$ as in (\ref{eqbern}). Besides, $\Gamma_\xi=(1/p) I_p$ and therefore
$\phi= \Omega \xi$ where $\Omega$ is given by~(\ref{eqomega}) with $\beta=p^{-1/2}$.
\end{ex}

\subsection{Orthogonal bases}
\label{expart2}

\subsubsection{Trigonometric basis}
\label{trigo}

The trigonometric family is defined by $\phi_0(x)=1$, $\phi_{2j-1}(x)=\sqrt{2}\sin(2\pi j x)$
and $\phi_{2j}(x)=\sqrt{2}\cos(2\pi j x)$ for all $j\geq 1$ and $x\in[0,1]$.
It is orthonormal with respect to the usual scalar product on $L_2([0,1])$.
Let $\theta\geq 0$ and consider the $(2p+1)\times (2p+1)$ matrix $A=$diag$\{1,\theta,\theta,\dots,\theta\}$.
One has
$$
c_{p,\theta}(u,v):=\;^t\!\phi(u)A\phi(v)= 1-\theta + \theta D_p(u-v)
$$
where $D_p$ is the Dirichlet kernel given by
$$
D_p(t)=\frac{\sin((2p+1)\pi t)}{\sin(\pi t)}.
$$
It is then clear that $c_{p,\theta}\in {\mathcal C}_\phi$ if $\theta\leq 1/(1-D_p^-)$
where $D_p^-:=\min_t D_p(t)$. Since $D_p^-<0$, it follows that that $\theta$ is upper
bounded by $1/(1-D_p^-)<1$ and thus $I_p\notin {\mathcal A}_\phi$ when $\phi$ is the trigonometric
family.
Besides, Proposition~\ref{proprho} shows that
the associated Spearman's rho is 
$$
\rho_{p,\theta} =\frac{6\theta}{\pi^2}\sum_{j=1}^p \frac{1}{j^2}.
$$
Numerical computations show that the maximum value is obtained for $p=2$ 
for which $D_2^-=-1$ and $\rho_{2,1/2}=15/(4\pi^2)$.
This bound can be increased by introducing
$$
\bar{c}_{q,\theta}(u,v):= \frac{1}{q} \sum_{p=0}^{q-1} c_{p,\theta}(u,v) = 1-\theta + \theta F_q(u-v)
$$
where $F_q$ is the Fej\'er kernel~\cite{Bases} defined as  
$$
F_q(x)=\frac{1}{q}\sum_{p=0}^{q-1} D_p(x) = \frac{1}{q}\left(\frac{\sin(qx/2)}{\sin(x/2)}\right)^2.
$$
Since this kernel is positive, it is readily seen that $\bar{c}_{q,\theta}$ is a density of copula
for all $\theta\in[0,1]$.
Besides, Lemma~\ref{lemagg} entails that $\bar{c}_{p,\theta}\in {\mathcal C}_\phi$ for all $\theta\in[0,1]$.
The associated Spearman's rho is 
$$
\bar{\rho}_{q,\theta} = \frac{1}{q} \sum_{p=0}^{q-1} \rho_{p,\theta} =\frac{6\theta}{\pi^2}\frac{1}{q} \sum_{p=0}^{q-1}\sum_{j=1}^p  \frac{1}{j^2} =
\frac{6\theta}{\pi^2}\left( \sum_{j=1}^{q-1} \frac{1}{j^2} - \frac{1}{q} \sum_{j=1}^{q-1} \frac{1}{j}
\right).
$$
Let us remark that $\bar{\rho}_{q,\theta} \to\theta$ as $q\to\infty$ and thus, arbitrary large
dependences can be modeled.

\subsubsection{The Haar basis}

For each positive integer $i$, let us denote by $J_i$ the interval
$J_i=\left[\frac{p_i}{ 2^{q_i-1}},\frac{p_i+1}{ 2^{q_i-1}}\right)$ where
$p_i$ and $q_i$
are the integers uniquely determined by
$i= 2^{q_i-1}+p_i$ and $0\leq p_i < 2^{q_i-1}$.
The Haar basis~\cite{Haar} is defined by
$$
\phi_0(t)=\mathbb{I}\{t\in[0,1]\}, \qquad
\phi_i(t)=2^{\frac{q_i-1}{2}}\left(\mathbb{I}\{t\in J_{2i}\}-\mathbb{I}\{t\in J_{2i+1}\}\right),
$$
for $i=1,\dots,p$. In the following, it is assumed that $p$ is a power of 2.
Let $\theta\geq 0$ and consider the $p\times p$ matrix $A=$diag$\{1,\theta,\theta,\dots,\theta\}$.
One has
\begin{equation}
\label{copuhaar}
 c_{p,\theta}(u,v):=\;^t\!\phi(u)A\phi(v)= 1-\theta + \theta K_p(u,v)
\end{equation}
where $K_p$ is the Dirichlet kernel associated to the Haar basis
$$
K_p(u,v)= p \sum_{i=1}^p \mathbb{I}\{(u,v) \in I_i^2\},
$$
and recall that $\{I_i,\; i=1,\dots,p\}$ is the equidistant partition of $[0,1]$ into $p$ intervals.
It is thus clear that $K_p$ can be rewritten as
$$
K_p(u,v)= p \sum_{i=1}^p \sum_{j=1}^p M_{i,j} \xi_i(u)\xi_j(v)
$$
where $M=I_p$ the $p\times p$ identity matrix and $\xi_i(u)=\mathbb{I}\{u\in I_i\}$.
It appears that $K_p$ is the density of a copula based on partition of unity.
In view of Proposition~\ref{propparti}, it follows that $K_p\in {\mathcal C}_\phi$. 
Finally, the Haar copula~(\ref{copuhaar}) can be interpreted as a linear mixture
of $K_p$ and the independent copula, leading to $c_{p,\theta}\in{\mathcal C}_\phi$
for all $\theta\in[0,1]$.
Thus, $I_p\in {\mathcal A}_\phi$ when $\phi$ is the Haar family.
Straightforward calculations show that the 
associated Spearman's rho is given by
$$
\rho_p(\theta)=\theta\left( 1-\frac{1}{p^2}\right).
$$
Let us remark that ${\rho}_p(\theta) \to\theta$ as $p\to\infty$ and thus, arbitrary large
dependences can be modeled.

\section{Conclusion and further work}
\label{seconcl}

We proposed a new family of copulas defined from a matrix and a family of orthogonal functions.
High dependences can be modeled without introducing singular components. It has also been shown that
this family can be used for approximating any density of copula. As a consequence, it appears
as a good tool for modeling bivariate data. 
To this aim, consider $(U_1,V_1),\dots,(U_n,V_n)$ independent copies of a random pair $(U,V)$ from
a density $c\in {\mathcal C}_\phi$ associated to a matrix $A$.
Assume that the function $\phi$ is known. Then, estimating $c$ reduces to estimating $A$.
Proposition~\ref{propT} provides two interpretations of $A$ in terms of a covariance matrices
and thus two possible estimators:
\begin{eqnarray*}
\hat A_{1,n} &=& \frac{1}{n} \sum_{i=1}^n \phi(U_i)\;^t\! \phi(V_i),\\
\hat A_{2,n} &=& \frac{1}{n} \sum_{i=1}^n (\phi(U_i)-e_1)\;^t\!(\phi(V_i)-e_1) + e_1 \;^t\!e_1.
\end{eqnarray*}
In both cases, the corresponding estimated density is given by 
$\hat c_{j,n}(u,v) = \;^t\!\phi(u) \hat A_{j,n} \phi(v)$, $j\in\{1,2\}$ and
can be simplified as
\begin{eqnarray*}
\hat c_{1,n}(u,v)& = & \frac{1}{n} \sum_{i=1}^n q(u,U_i) q(v,V_i),\\
\hat c_{2,n}(u,v)&= &1+ \frac{1}{n} \sum_{i=1}^n (q(u,U_i)-1) (q(v,V_i)-1).
\end{eqnarray*}
From the theoretical point of view, both estimators are unbiased and asymptotically Gaussian.
Our future work will consist in investigating their finite sample behavior 
compared to rank-based methods~\cite{genest}.
% 
% \section{Conclusion}
% %---------------------------------------------------------------%
% 
% \begin{itemize}
% \item We proposed a new family of copulas defined from matrices,
% \item The family includes some known families (FGM, cubic sections,...),
% \item The family contains high dependence copulas (Haar),
% \item Other orthonormal families should be studied (polynomial...),
% \item Dependence properties have to be more deeply studied.
% \end{itemize}

%**********************************************************************
%**********************************************************************

\section{Appendix}
 %---------------------------------------------------------------%
 
 \subsection{Proofs of main results}

\paragraph{Proof of Proposition~\ref{propdens}.}

Let $c\in {\mathcal C}_\phi$. Clearly, $c(u,v)\geq 0$
for all $(u,v)\in [0,1]^2$ from the definition of ${\mathcal A}_\phi$.
It only remains to prove that the margins of $c$ are standard uniform distributions.
To this end, let us remark that
\begin{equation}
 \label{intphi}
 \int_0^1 \phi(v)dv = e_1,
\end{equation}
since for all $j=1,\dots,p$, we have
$$
 \int_0^1 \phi_j(v)dv = \langle e_1,e_j \rangle = \delta_{1j},
$$
where $\delta_{ij}=1$ if $i=j$ and $\delta_{ij}=0$ otherwise. As a consequence,
$$
\int_0^1 c(u,v)dv = \;^t\phi(u) A \int_0^1 \phi(v)dv = \;^t\phi(u) A e_1 = \;^t\phi(u) e_1 = \phi_1(u)=1.
$$
The proof of $\int_0^1 c(u,v)du =1$ is similar. \hfill $\Box$\\

\paragraph{Proof of Proposition~\ref{propT}.}

Let $c\in {\mathcal C}_\phi$ such that 
$c(u,v)=\;^t\!\phi(u) A \phi(v)=\;^t\!\phi(u) B \phi(v)$
for all $(u,v)\in [0,1]^2$ and for some $(A,B)\in {\mathcal A}_\phi^2$.
It follows that 
$$
c(u,v) \phi(u)\;^t\!\phi(v)=\phi(u)\;^t\!\phi(u)A\phi(v)\;^t \! \phi(v)= \phi(u)\; ^t\!\phi(u) B\phi(v)\; ^t \!\phi(v),
$$
and integrating with respect to $u$ and $v$ yields $T_\phi(c)=A=B$ in view of~(\ref{ortho}). 
Remark also that, if moreover, $c\in {\mathcal C}_\phi$ with associated matrix $A\in {\mathcal A}_\phi$ then 
$$
{\mathbb E}_c(\phi(U))= \int_0^1\int_0^1 \phi(u)\;^t\!\phi(u)A\phi(v) dudv = A \int_0^1\phi(v)dv = A e_1=e_1
$$
in view of (\ref{ortho}) and (\ref{intphi}). As a consequence, the matrix $A$ can also be interpreted as
$A=$cov$(\phi(U),\phi(V))+ e_1\; ^t\! e_1$.\hfill $\Box$\\

\paragraph{Proof of Proposition~\ref{propA}.}

It is clear that ${\mathcal A}_\phi$ is convex. Let us prove that
$({\mathcal A}_\phi,\times)$ is a semi-group. Since the product $\times$ is associative,
it only remains to establish that $(A,B)\in {\mathcal A}_\phi^2$ entails $A\times B\in {\mathcal A}_\phi$.
First, (\ref{ortho}) entails
\begin{eqnarray*}
\;^t\phi(u) AB  \phi(v)&=& \;^t\phi(u) A \left\{\int_0^1 \phi(y)\;^t\phi(y)dy\right\}  B  \phi(v)\\
&=& \int_0^1 \left\{ \;^t\phi(u) A \phi(y)\right\} \left\{\;^t\phi(y) B \phi(v)\right\}dy
\geq  0,
\end{eqnarray*}
and second it is easily seen that $AB e_1=e_1$ and $^t\!(AB)e_1=e_1$.
Finally, $I_p \in {\mathcal A}_\phi$ and $({\mathcal A}_\phi,\times)$ is a semi-group
both imply that $({\mathcal A}_\phi,\times)$ is a monoid. \hfill $\Box$\\

\paragraph{Proof of Proposition~\ref{propC}.}

It is clear that ${\mathcal C}_\phi$ is convex. Let us prove that
$({\mathcal C}_\phi,\star)$ is a semi-group. Since the product $\star$ is associative
(see~\cite{Darsow}, Theorem~2.4),
it only remains to establish that $(c_A,c_B) \in {\mathcal C}_\phi^2$ entails
$c_A\star c_B \in {\mathcal C}_\phi$:
\begin{eqnarray*}
c_A\star c_B(u,v)
&=&\int_0^1 \; ^t\!\phi(u) A \phi(s)\; ^t\! \phi(s) B \phi(v)ds,\\
&=& ^t\!\phi(u) A \left\{\int_0^1 \phi(s)\;^t\!\phi(s) ds\right\} B\phi(v)
= ^t\!\phi(u) AB \phi(v).
\end{eqnarray*}
This proves that $(c_A,c_B) \in {\mathcal C}_\phi^2$ and is associated to the matrix $AB\in {\mathcal C}_\phi$, see Proposition~\ref{propA}. 
Finally, $q \in {\mathcal C}_\phi$ and $({\mathcal C}_\phi,\star)$ is a semi-group
both imply that $({\mathcal C}_\phi,\star)$ is a monoid. \hfill $\Box$\\

\paragraph{Proof of Proposition~\ref{proprho}.}

Let $C$ denote the copula associated to the density
$c\in{\mathcal C}_\phi$. Introducing $\Psi(u)=\int_0^u \phi(t)dt$
and $\gamma=\int_0^1 \Psi(u)du$, we have
$C(u,v)=\;^t\!\Psi(u)A\Psi(v)$ and~(\ref{eqrho}) leads to $\rho=12\;^t\!\gamma A \gamma -3$.
A partial integration yields $\gamma=e_1-\mu$ and 
$$
\;^t\!(e_1-\mu) A (e_1-\mu) = \;^t\!e_1 A e_1 - \;^t\!\mu(A+ \;^t\!A) e_1 + \;^t\!\mu A \mu
=1 - 2 \;^t\!\mu_1 + \;^t\!\mu A \mu
$$
with $\mu_1=\int_0^1 t dt=1/2$ and the result follows. \hfill $\Box$\\

\paragraph{Proof of Proposition~\ref{proptau}.}

Let $C$ denote the copula associated to the density
$c\in{\mathcal C}_\phi$. Standard algebra gives
\begin{eqnarray*}
 \int_0^1 \!\!\! \int_0^1 C(u,v) dC(u,v) 
 & = & \int_0^1 \!\!\! \int_0^1 \;^t\!\Psi(u)A \Psi(v) \;^t\!\phi(u)A \phi(v)dudv \\
 & = & \int_0^1 \!\!\! \int_0^1 \;^t\!\Psi(v)\;^t\!A \Psi(u) \;^t\!\phi(u)A \phi(v)dudv \\
 & = & \int_0^1 \;^t\!\Psi(v)\;^t\!A \Theta A \phi(v) dv \\
 & = & \int_0^1 \mbox{tr} \left( \;^t\!\Psi(v)\;^t\!A \Theta A \phi(v)\right) dv \\
 & = & \int_0^1 \mbox{tr} \left( \;^t\!A \Theta A \phi(v)\;^t\!\Psi(v)\right) dv \\
 & = & \mbox{tr} \left( \;^t\!A \Theta A \;^t\!\Theta\right).
\end{eqnarray*}
Besides, a partial integration shows that $\;^t\!\Theta=e_1 \;^t\!e_1 -\Theta$.
It follows that
\begin{eqnarray*}
\mbox{tr} \left( \;^t\!A \Theta A \;^t\!\Theta\right) &=& 
\mbox{tr} \left( \;^t\!A \Theta A e_1 \;^t\!e_1\right) - \mbox{tr} \left( \;^t\!A \Theta A \Theta\right)\\
&=& 
\mbox{tr} \left( \;^t\!A \Theta e_1 \;^t\!e_1\right) - \mbox{tr} \left( \;^t\!A \Theta A \Theta\right)\\
&=& 
\mbox{tr} \left(  \Theta e_1 \;^t\!e_1\;^t\!A\right) - \mbox{tr} \left( \;^t\!A \Theta A \Theta\right)\\
&=& 
\mbox{tr} \left( \Theta e_1 \;^t\!e_1\right) - \mbox{tr} \left( \;^t\!A \Theta A \Theta\right)\\
&=& 
\Theta_{1,1}  - \mbox{tr} \left( \;^t\!A \Theta A \Theta\right)\\
&=& 
1/2  - \mbox{tr} \left( \;^t\!A \Theta A \Theta\right).
\end{eqnarray*}
The conclusion follows from~(\ref{eqtau}). \hfill $\Box$\\

\paragraph{Proof of Proposition~\ref{proplambda}.}

Let $C$ denote the copula associated to the density
$c\in{\mathcal C}_\phi$. Recalling that $\Psi(u)=\int_0^u \phi(t)dt$, we have
$C(u,v)=\;^t\!\Psi(u)A\Psi(v)$. In view of~(\ref{intphi}),
$\Psi(1)=e_1$ leading to $C(1,1)=1$ and thus~(\ref{eqlambda}) can be rewritten as
$$
\lambda=2- \lim_{u\rightarrow 1^-}\frac{C(u,u)-C(1,1)}{u-1} = 2 - \left.\frac{\partial C(u,u)}{\partial u}\right|_{u=1}.
$$
Straightforward calculations show that
\begin{eqnarray*}
 \left.\frac{\partial C(u,u)}{\partial u}\right|_{u=1}
=  \;^t\!\phi(1) (A + \;^t\! A) \Psi(1) 
=  \;^t\!\phi(1) (A + \;^t\! A) e_1 
=  2 \;^t\!\phi(1) e_1 =2,
\end{eqnarray*}
and consequently $\lambda=0$. \hfill $\Box$\\

\paragraph{Proof of Proposition~\ref{propP}.}

Let us suppose first that $I_p\in {\mathcal A}_\phi$.
Let us first remark that, from Proposition~\ref{propC}, ${\mathcal C}_\phi$ is a convex and closed
subset of $L_2([0,1]^2)$. Second, it is clear from Lemma~\ref{lemP}(i) that $P_\phi$ is idempotent:
$$
P_\phi(P_\phi(c)) = q \star P(c)\star q = q \star q \star c\star q \star q = q \star c\star q = P(c),
$$
for all $c\in L_2([0,1]^2)$ since $q\star q=q$ from (\ref{ortho}).
Third, let $c\in L_2([0,1]^2)$ and $s\in {\mathcal C}_\phi$
with associated matrix $A$. Our aim is to prove that $\prec c-P_\phi(c),s\succ=0$.
In view of Lemma~\ref{lemps}, 
$$
\prec c-P_\phi(c),s\succ = \mbox{tr}(T_\phi(c)A) - \mbox{tr}(T_\phi(c)A) =0
$$
and the direct part of the result is proved.
Conversely, if $P_\phi$ is a projection on ${\mathcal C}_\phi$ then, necessarily, $P_\phi(q)\in {\mathcal C}_\phi$. Besides, $P_\phi(q)=q$ and thus $q\in {\mathcal C}_\phi$ entailing $I_p\in {\mathcal A}_\phi$.
The converse part of the result is proved.
\hfill $\Box$

\paragraph{Proof of Corollary~\ref{coroprop}.}

First, it is clear that $\|c-P_\phi(c)\|\to 0$ as $p\to\infty$, since, in view of Theorem~\ref{propproj}(i),
$\tilde{c}_p := P_\phi(c)$ can be interpreted as a $L_2-$ projection of $c$. 
From~(\ref{eqrho}), it follows that
$$
\rho(\tilde{c}_p) -\rho(c) = \int_{[0,1]^4} (\tilde{c}_p(x,y)-c(x,y)) \mathbb{I}\{x\leq u\} \mathbb{I}\{y\leq v\} dxdydudv
$$
and Cauchy-Schwarz inequality yields
$$
|\rho(\tilde{c}_p) -\rho(c)|\leq \|\tilde{c}_p -c \|  \left( \int_{[0,1]^4}  \mathbb{I}\{x\leq u\} \mathbb{I}\{y\leq v\} dxdydudv\right)^{1/2} = \frac{1}{2} \|\tilde{c}_p -c \|
$$
and the conclusion follows. \hfill $\Box$
% 
% \paragraph{Proof of Proposition~\ref{propNA}.}
% The proof is a consequence of some simple calculations
% ${\mathbb E}(q(u,U))= 1$, ${\mathbb E}(q(u,U)q(v,V))=c(u,v)$ and on the Central-Limit Theorem.\hfill $\Box$
% 
% \paragraph{Proof of Proposition~\ref{propnorm}.}
% The proposition is the consequence of a classical result of the $L_2-$ projection:
% $$
% \| \hat c_{j,n} - c \|^2 = \|\hat c_{j,n} -\tilde{c}_p \|_2 + \| \tilde{c}_p -c \|^2.
% $$
% The result follows from the remark that
% $$
% \|\hat c_{j,n} -\tilde{c}_p \|_2 = \|\hat A_{j,n} - A_p\|_F^2. 
% $$
% \hfill $\Box$

\paragraph{Proof of Proposition~\ref{propparti}.}

 Let 
$\psi=H\xi$ and $B=p \;^t\!H^{-1} M H^{-1}$.
First, it is easily seen that 
$$
c(u,v)= p \;^t\!\xi(u) M \xi(v) = \;^t\!\psi(u) B \psi(v),
$$
with 
$$
\psi_1(t)= \sum_{i=1}^p \xi_i(t) = 1
$$
and
$$
\int_0^1 \psi(t)dt= H \int_0^1 \xi(t)dt = \frac{1}{p} H s = e_1. 
$$
Besides, since $\xi_i\geq 0$ for all $i=1,\dots,p$ and $M$ is a doubly stochastic matrix,
it is clear that $c(u,v)\geq 0$ for all $(u,v)\in[0,1]^2$.
Second, standard algebra shows that $H^{-1}e_1= s/p$ and $\;^t\!H^{-1}s=e_1$. As a consequence,
$$
Be_1 = p \;^t\!H^{-1} M H^{-1}e_1 = \;^t\!H^{-1} M s =  \;^t\!H^{-1} s =e_1
$$
and similarly $\;^t\!Be_1=e_1$ leading to $B\in {\mathcal A}_\psi$.
Lemma~\ref{lemortho} entails that $A:=\Gamma^{1/2} B \Gamma^{1/2} \in {\mathcal A}_\phi$
where $\phi:=\Gamma^{-1/2} \psi$ fullfills the conditions of Definition~\ref{def1}
and the density of copula
$$
c(u,v)= p \;^t\!\xi(u) M \xi(v) = \;^t\!\psi(u) B \psi(v)= ^t\!\phi(u)A\phi(v)
$$ 
belongs to ${\mathcal C}_\phi$.\hfill $\Box$

 \subsection{Proofs of auxiliary results}
\label{proofaux}

\paragraph{Proof of Lemma~\ref{lemortho}.}

Let us first remark that $\Gamma e_1=e_1$ and $^t \Gamma e_1=e_1$.
Consequently, there exists a square root $\Gamma^{1/2}$ of $\Gamma$ such that
$\Gamma^{1/2} e_1=e_1$ and $^t \Gamma^{1/2} e_1=e_1$. It follows that
$Ae_1=e_1$,  $\;^t\!Ae_1=e_1$
and that, for all $(u,v)\in[0,1]^2$, $^t\!\phi(u)A\phi(v)=\;^t\!\psi(u)B\psi(v)$, it is thus clear that $A\in {\mathcal A}_\phi$
with $\phi_1(x)=1$ for all $x\in[0,1]$ and $\phi$ is orthonormal.\hfill $\Box$\\

\paragraph{Proof of Lemma~\ref{lemagg}.}
Let $c_p\in {\mathcal C}_\phi$ with associated matrix $A$ and consider the 
density of copula defined for all $(u,v)\in[0,1]^2$ and $q\geq 1$ by
\begin{eqnarray*}
 \bar{c}_{q}(u,v)&:=& \frac{1}{q} \sum_{p=1}^{q} c_{p}(u,v) \\
 &=&\frac{1}{q} \sum_{p=1}^{q} \sum_{i=1}^{q}\sum_{j=1}^{q} A_{i,j} \phi_i(u)\phi_j(v) \mathbb{I}\{i\leq p\}\mathbb{I}\{j\leq p\}\\
 &=&  \sum_{i=1}^{q}\sum_{j=1}^{q} A_{i,j} \left(\frac{1}{q} \sum_{p=1}^{q} \mathbb{I}\{\max(i,j) \leq p\}\right)\phi_i(u)\phi_j(v)\\
 &=&  \sum_{i=1}^{q}\sum_{j=1}^{q} A_{i,j} \left( \frac{q+1-\max(i,j)}{q} \right)\phi_i(u)\phi_j(v)\\
  &=:&  \sum_{i=1}^{q}\sum_{j=1}^{q} B_{i,j} \phi_i(u)\phi_j(v).
\end{eqnarray*}
It is clear from its definition that $\bar{c}_{q}$ is a density of copula. Therefore, $\bar{c}_{q}\in {\mathcal C}_\phi$ and the result is proved.\hfill $\Box$\\

\paragraph{Proof of Lemma~\ref{lemP}.}

(i) Let $c\in L_2([0,1]^2)$ and suppose $I_p\in {\mathcal A}_\phi$. Then
$q \in {\mathcal C}_\phi$ and
\begin{eqnarray*}
\;^t\!\phi(u) T_\phi(c)\phi(v) &=& \;^t\!\phi(u) \left\{\int_0^1\!\!\! \int_0^1  \phi(x)c(x,y)\;^t\!\phi(y)dxdy\right\}\phi(v)\\
&=& \int_0^1\!\!\! \int_0^1 q(u,x) c(x,y) q(y,v) dxdy\\
&=&(q \star c\star q)(u,v),
\end{eqnarray*}
which proves that $\;^t\!\phi(u) T_\phi(c)\phi(v)\geq 0$ for all $(u,v)\in[0,1]^2$.
Moreover,
\begin{eqnarray*}
T_\phi(c) e_1&=& \int_0^1\!\!\! \int_0^1 \phi(x)c(x,y)\;^t\! \left\{\phi(y) e_1\right\} dydx
= \int_0^1\!\!\! \int_0^1 \phi(x)c(x,y) dydx\\
&=&\int_0^1 \phi(x)dx
=e_1,
\end{eqnarray*}
in view of (\ref{intphi}) and similarly $\;^t\!T_\phi(c) e_1=e_1$.
As a conclusion, $T_\phi(c)\in {\mathcal A}_\phi$ and thus $P_\phi(c)\in {\mathcal C}_\phi$.\\
(ii) Conversely, if $P_\phi(c)\in {\mathcal C}_\phi$ for all $c\in L_2([0,1]^2)$, we have in particular
$P_\phi(q)=q\in {\mathcal C}_\phi$ and thus $I_p\in{\mathcal A}_\phi$.
\hfill $\Box$

\paragraph{Proof of Lemma~\ref{lemps}.}

(i) Let $c_2\in L_2([0,1]^2)$ and let $c_1\in {\mathcal C}_\phi$ with associated matrix $A$. By definition,
\begin{eqnarray*}
\prec c_1,c_2\succ &=& \int_0^1\!\!\!\int_0^1 c_1(u,v) c_2(v,u) du dv\\
&=&  \int_0^1\!\!\!\int_0^1 \;^t\!\phi(u) A \phi(v)c_2(v,u)  du dv \\
&=& \int_0^1\!\!\!\int_0^1 \mbox{tr} \left\{ \;^t\!\phi(u) A \phi(v)c_2(v,u)\right\} du dv \\
&=& \int_0^1 \!\!\!\int_0^1 \mbox{tr} \left\{  A \phi(v)c_2(v,u)\;^t\!\phi(u) \right\} du dv 
= \mbox{tr}\left\{ A T_\varphi(c_2)  \right\}.
\end{eqnarray*}
(ii) If, moreover, $c_2\in {\mathcal C}_\phi$ with associated matrix $B$ then $T_\varphi(c_2)=B$
and the conclusion follows. \hfill $\Box$


\begin{thebibliography}{99}
%**********************************************************************
%**********************************************************************



  
\bibitem{Nous} Amblard, C. and Girard, S., 2002.
Symmetry and dependence properties within a semiparametric family
of bivariate copulas. {\em Journal of Nonparametric Statistics}, {\bf 14},
715--727.

\bibitem{Nous2} Amblard, C. and Girard, S., 2005.
Estimation procedures for a semiparametric family of bivariate copulas.
{\em Journal of Computational and Graphical Statistics}, {\bf 14}, 1--15.

\bibitem{Metrika} Amblard, C. and Girard, S., 2009.
A new extension of bivariate FGM copulas.
{\em Metrika}, {\bf 70}, 1--17.

\bibitem{Darsow} Darsow, W., Nguyen, B. and Olsen, E., 1992.
Copulas and Markov processes.
{\em Illinois Journal of Mathematics}, {\bf 36}, 600--642.

\bibitem{DUR} Durante, F., 2006.
A new class of symmetric bivariate copulas. {\em Journal of Nonparametric Statistics}, {\bf 18}, 499--510.

\bibitem{DUR2} Durante, F., Koles{\'a}rov{\'a}, A., Mesiar, R. and Sempi, C., 2008.
Semilinear copulas.
{\it Fuzzy Sets and Systems}, {\bf 159}, {63--76}.
  
\bibitem{GRO} Durrleman, V., Nikeghbali, A. and Roncalli, T., 2000.
Copulas approximation and new families.
{\em Technical report.}
{\tt http://dx.doi.org/10.2139/ssrn.1032547}.

\bibitem{FAR} {Farlie, D. G. J., 1960. }
The performance of some correlation coefficients for a general bivariate distribution.
{\em Biometrika}, {\bf 47}, 307--323.

\bibitem{FISC} {Fischer, M. and Klein, I., 2007.}
Constructing symmetric generalized FGM copulas by means
of certain univariate distributions,
{\em Metrika}, {\bf 65}, 243--260.

\bibitem{genest} Genest, C. and Favre, A.-C., 2007.
{Everything you always wanted to know about copula modeling but were afraid to ask},
{\em Journal of Hydrologic Engineering},
{\bf 12}, {347--368}.

\bibitem{JNS10} Genest, C., Neslehov\`a, J. and Ghorbal, N. B., 2010.
Spearman's footrule and Gini's gamma: a review with complements,
{\it Journal of Nonparametric Statistics}, {\bf 22}, 937--954.

\bibitem{GUMBEL} {Gumbel, E. J., 1960. }
Bivariate exponential distributions.
{\em Journal of the American Statistical Association}, {\bf 55}, 698--707.

\bibitem{Bases} {Hoffman, K., 1988.}
{\it Banach spaces of analytic functions}, Dover.


 \bibitem{IFGM2} {Huang, J. S. and Kotz, S., 1984. }
 Correlation structure in iterated Farlie-Gumbel-Morgenstern distributions,
 {\it Biometrika},
 {\bf 71}, {633--636}.
 
 \bibitem{HUE} {Huang, J. S. and Kotz, S., 1999. }
 Modifications of the Farlie-Gumbel-Morgenstern distribution.
 A tough hill to climb.
 {\em Metrika}, {\bf 49}, 135--145.
 
\bibitem{JOE} {Joe, H., 1997. }
{\em Multivariate models and dependence concepts}.
Monographs on statistics and applied probability, {\bf 73},
Chapman \& Hall.

\bibitem{IFGM1} Johnson, N. and Kotz, S., 1975.
On some generalized Farlie-Gumbel-Morgenstern distributions,
{\it Communications in Statistics-Theory and Methods},
{\bf 4}, {415--427}.

 \bibitem{LAI} {Lai, C. D. and Xie, M., 2000. }
A new family of positive quadrant dependence bivariate distributions.
{\em Statistics and Probability Letters}, {\bf 46}, 359--364.

\bibitem{LMST} Li X., Mikusi\'nski, P., Sherwood, H. and Taylor, M. D., 1997.
On approximation of copulas. In: Bene$\check{s}$ V, $\check{S}$t$\check{e}$p\`an J (eds)
{\it Distributions with given marginals and moment problems}. Kluwer Academic Publishers, Dordrecht.

\bibitem{LMT} Li X., Mikusi\'nski, P. and Taylor, M. D., 1998.
Strong approximation of copulas. 
{\em Journal of Mathematical Analysis and Applications}, {\bf 225}, 608--623.

\bibitem{IFGM3} Lin, G. D., 1987.
 Relationships between two extensions of Farlie-Gumbel-Morgenstern distribution,
 {\it Annals of the Institute of Statistical Mathematics},
 {\bf 39}, {129--140}.
  
\bibitem{MORG} {Morgenstern, D., 1956. }
Einfache Beispiele zweidimensionaler Verteilungen.  
{\em Mitteilungsblatt f\"ur Mathematische Statistik}, {\bf 8}, 234--235.

\bibitem{NELSEN97} {Nelsen, R. B., Quesada-Molina, J. J. and Rodr\'{\i}guez-Lallena, J. A., 1997. } Bivariate copulas with cubic sections.
{\em Journal of Nonparametric Statistics}, {\bf 7}, 205--220.

\bibitem{JNS98} Nelsen, R. B., 1998.
Concordance and Gini's measure of association, {\it Journal of Nonparametric Statistics}, {\bf 9},
227--238.

\bibitem{NELSEN06} {Nelsen, R. B., 2006. }
{\em An introduction to copulas, Second Edition}.  
Springer Series in Statistics, Springer.


\bibitem{LLALENA} {Rodr\'{\i}guez-Lallena, J. A., 1992.} 
{\em Estudio de la compabilidad y dise$\tilde{n}$o de nuevas familias
en la teoria de c\'opulas. Aplicaciones}.
Tesis doctoral, Universidad de Granada.

\bibitem{LLALENA2} {Rodr\'{\i}guez-Lallena, J. A. and \'Ubeda-Flores, M., 2004. }
A new class of bivariate copulas.
{\em Statistics and Probability Letters}, {\bf 9}, 315--325.

\bibitem{QUESADA} {Quesada-Molina, J. J. and Rodr\'{\i}guez-Lallena, J. A., 1995. }
Bivariate copulas with quadratic sections.
{\em Journal of Nonparametric Statistics}, {\bf 5}, 323--337.

\bibitem{Haar} Ruch, D. K., Van Fleet, P. J., 2009.
{\it Wavelet theory: An elementary approach with applications}, John Wiley and Sons.

\bibitem{sancetta}   {Sancetta, A. and Satchell, S., 2004.}
The Bernstein copula and its applications to modeling and approximations of multivariate distributions.
{\em Econometric Theory}, {\bf 20}, {535--562},

\bibitem{SIGMA} Schweizer, B. and Wolff, E. F, 1981. 
On nonparametric measures of dependence for random variables.
{\em The Annals of Statistics}, {\bf 9}, 879--885.

\bibitem{SKLAR} {Sklar, A., 1959. }
Fonctions de r\'epartition \`a $n$ dimensions et leurs marges.
{\em Publ. Inst. Statist. Univ. Paris}, {\bf 8}, 229--231.

\end{thebibliography}
\end{document}